\documentclass[twocolumn]{autart}
\usepackage{graphicx}
\usepackage{amsmath}
\usepackage{amssymb}
\usepackage{float}

\def\e{{\mathbf{e}}}          
\def\f{{\mathbf{f}}}          
          
\def\h{{\mathbf{h}}}

\def\v{{\mathbf{v}}}          
\def\w{{\mathbf{w}}}          
\def\x{{\mathbf{x}}}          
\def\y{{\mathbf{y}}}

\def\Dpartial#1#2{ \frac{\partial #1}{\partial #2} }

\DeclareMathSymbol{\Complex}{\mathalpha}{AMSb}{"43}
\DeclareMathSymbol{\Real}{\mathalpha}{AMSb}{"52}

\newcommand{\pderiv}[2]{\frac{\partial #1}{\partial #2}}
\newcommand{\pderivsq}[2]{\frac{\partial^{2} #1}{\partial #2^{2}}}

\begin{document}
\begin{frontmatter}
\title{Efficient grid-based Bayesian estimation of nonlinear low-dimensional systems with sparse non-Gaussian PDFs}
\author[UCSD]{Thomas R.~Bewley}\ead{bewley@ucsd.edu}
\author[Sheffield]{Atul S.~Sharma}\ead{a.s.sharma@sheffield.ac.uk}
\address[UCSD]{Flow Control \& Coordinated Robotics Labs, Dept of MAE, UC San Diego, La Jolla, CA 92093-0411, USA}
\address[Sheffield]{Department of Automatic Control and Systems Engineering, University of Sheffield, S1 3JD, UK}


\begin{abstract} 
Bayesian estimation strategies represent the most fundamental formulation of the state estimation problem available, and apply readily to nonlinear systems with non-Gaussian uncertainties.
The present paper introduces a novel method for implementing grid-based Bayesian estimation which largely sidesteps the severe computational expense that has prevented the widespread use of such methods.
The method represents the evolution of the probability density function (PDF) in phase space, $p_{\x}(\x',t)$, discretized on a fixed Cartesian grid over {\it all} of phase space, and consists of two main steps:
(i) Between measurement times, $p_{\x}(\x',t)$ is evolved via numerical discretization of the Kolmogorov forward equation, using a Godunov method with second-order corner transport upwind correction and a total variation diminishing flux limiter;
(ii) at measurement times, $p_{\x}(\x',t)$ is updated via Bayes' theorem.
Computational economy is achieved by exploiting the localised nature of $p_{\x}(\x',t)$. An ordered list of cells with non-negligible probability, as well as their immediate neighbours, is created and updated, and the PDF evolution is tracked {\it only} on these active cells.
\end{abstract}
\end{frontmatter}

\section{Introduction}
Bayesian estimation strategies are the most general class of solutions to the state estimation problem, and apply readily to nonlinear systems where information about the state is represented by a probability density function (PDF) of general form.
In this paper we introduce a novel, computationally cheap method for implementing grid-based Bayesian estimation that exploits the fact that the PDF is usually negligible in most of phase space, while avoiding many of the disadvantages of other methods.
The idea of grid-based Bayesian estimation dates back at least to Stratonovich (1959,1960). The equations underlying the algorithm are laid out clearly in Jazwinski (1970, p.~164), and are summarised below. However, numerical implementation of these equations has only been attempted sporadically in the half century since, for instance by Kramer et al. (1988), Terwiesch \& Agarwal (1994) and Ungarala et al. (2006).
Grid-based Bayesian methods typically suffer from the twin burdens of high computational cost and a finite domain size; indeed, Arulampalam et al. (2002), in their otherwise insightful review of particle filter methods, all but dismiss grid-based methods in \S IV.B. We believe that this level of pessimism on this class of methods is unwarranted.
The algorithm developed, dubbed GBEES (Grid-based Bayesian Estimation Exploiting Sparsity), provides a means of efficient computation by building on an accurate integration scheme for hyperbolic systems, and a novel gridding scheme over \emph{all} of phase space.

\section{Grid-based Bayesian estimation exploiting sparsity}
\label{sec:nummethsparse}

Consider the state estimation of the nonlinear system
\begin{equation}
\frac{d\x}{dt} = \f(\x,\w),\qquad \y = \h(\x,\v).
\label{sysnon}
\end{equation}
The grid-based Bayesian estimation method is best visualized as an evolution of the PDF of the state estimate $\hat\x$ discretized on a fixed grid over all of phase space $\Real^n$; assuming the state $\x$ develops according to the nonlinear equation \eqref{sysnon}, the method consists of two relatively straightforward steps (for details, see Jazwinski 1970, p.~164):\\[0.05in]
\noindent (i) Between measurement times, the PDF itself, $p_{\x}(\x',t)$, is marched via discretization of the {\it Kolmogorov forward equation} (also called the {\it Fokker-Planck equation})
\begin{equation}
\Dpartial{p_{\x}(\x',t)}{t} = - \Dpartial{f_i(\x',t)\, p_{\x}(\x',t)}{x_i'} +
\frac{1}{2} \frac{\partial^2 q_{ij}\,p_{\x}(\x',t)}{\partial x_i' \,\partial x_j'}, \label{eq:PDEgovPDF}
\end{equation}
where summation over repeated indices is implied and $q_{ij}$ is the $(i,j)$th element of the spectral density, $Q$, of the state disturbances (note that, in the special case that $Q$ is diagonal and the state disturbance is independent of $\x$, this is just a diffusion term).
Risken (2002) reviews a number of methods for solving this equation, including analytic methods for special cases and eigenfunction expansions (focusing on the stationary solution), however state-space and time discretisation for a non-stationary solution (as performed here) is only mentioned briefly, and assumed to apply only to finite domains.
An accurate numerical method for marching this equation in time in the case that $Q=0$ is outlined in \S \ref{sec:nummeth}; adding an appropriate term to this discretization to apply diffusion to the PDF (to account for Gaussian state disturbances) is straightforward, as discussed in \S \ref{sec:nummethdiff}.\\[0.05in]
\noindent (ii) At the measurement times $t_k$, the PDF is updated via Bayes' theorem (Bayes, 1763),
\begin{equation}
p_{\x}(\x',t_{k^+})=\frac{p_{\y}(\y_k|\x')\,p_{\x}(\x',t_{k^-})}{C},
\label{eq:BayesforGBEES}
\end{equation}
where $p_{\x}(\x',t_{k^+})$ denotes the {\it a posteriori} PDF (after accounting for the measurement $\y_k$),
$p_{\y}(\y_k|\x')$ denotes the uncertainty associated with the measurement (which may or may not be Gaussian in $\x'$),
$p_{\x}(\x',t_{k^-})$ denotes the {\it a priori} PDF (before accounting for the measurement $\y_k$),
and $C$ is an appropriate normalization constant, which is selected for every measurement update to normalize
the discretization of $p_{\x}(\x',t_{k^+})$ such that its integral over phase space is unity.

To understand how the continuous state-space is discretized, recall that the {cumulative distribution function} ({CDF}) of a random real vector $\x$, denoted $f_\x(\underline\x)$, maps $\underline\x\in\Real^{n}$ to the real interval $[0,1]$ that monotonically increases in each of the components of $\underline\x$, and is defined
\begin{equation*}
f_\x(\underline\x)=P(x_1\le \underline x_1, x_2\le \underline x_2, \ldots, x_n\le \underline x_n),
\end{equation*}
where $\underline \x$ is some particular value of the random vector $\x$ and $P(S)$ denotes a {probability measure} that the conditions stated in $S$ are true. For any random vector $\x$ whose CDF is differentiable everywhere, the probability density function (PDF) $p_\x(\x')\ge 0$ is a scalar function of $\x'$ defined such that
\begin{align*}
&f_\x(\underline \x) = \int_{-\infty}^{\underline x_1}\int_{-\infty}^{\underline x_2}\cdots \int_{-\infty}^{\underline x_n}
p_\x(\x') \,dx'_1\,dx'_2\,\cdots dx'_n, \\
&\Leftrightarrow \qquad p_\x(\x') = \frac{\partial^n \! f_\x(\underline \x)}{\partial \underline x_1\,\partial \underline x_2 \cdots \partial \underline x_n}
\Big|_{\underline \x=\x'}.
\end{align*}
For small $|\Delta \x'|$, the quantity $p_{\x}(\x') \Delta x'_1 \,\Delta x'_2 \cdots \Delta x'_n$ represents the probability that the random vector $\x$ takes some value within a small rectangular region centered at the particular value $\x'$ and of width $\Delta x'_i$ in each coordinate direction $\e_i$.

The method we have developed maintains a list of active cells on the grid over all of phase space in order to limit both the computational effort and the memory storage required in the numerical simulation.
This list includes all cells in the discretization for which the PDF is greater than a given threshold, as well as all cells which, though they may or may not themselves exceed this threshold, are either one of the two immediate neighbor cells, in each of the $n$ coordinate directions, of those cells which exceed the threshold, or are one of the four neighbor cells, in each of the ${}_nC_2$ pairs of coordinate directions, which touch a corner of the cells which exceed the threshold.
At each time step, cells are added to and removed from this list as appropriate, and the fluxes initialized and updated on every interior boundary between adjacent cells in the list.
When performing a computation restricted to an evolving list of active grid cells of this sort, the relative position of the various cells in the list is needed frequently.
This may be determined efficiently by keeping in each list record a pointer to the two immediate neighbor cells in each coordinate direction in the list, if these neighbor cells are present in the list, or to a null record if not, and updating these pointers appropriately as records are added to and removed from the list\footnote{The four neighbor cells, in each pair of coordinate directions, which touch the corner of a given cell may be found by referencing the neighbor cell of a neighbor cell.}.  These pointers interconnecting the list facilitate rapid computation of the numerical discretization given in \eqref{eq:hyper_conD} in the next section.
The most expensive step in maintaining this list of neighbor cells is making the appropriate connections when a new record is added to the list. Though this may be accomplished by scanning the entire list, this approach becomes prohibitively expensive as the length of the list grows to thousands of cells. Instead, we keep the list ordered by its indices (e.g., in a phase space with $n=3$, ordered first by $i$, then by $j$, then by $k$), and store the elements of the list as a binary tree. This allows the time-limiting search step to proceed at ${\mathcal{O}}(N\log N)$ operations, where $N$ is the number of list elements.  Conveniently, this list ordering and searching can be handled using the C++ Standard Template Library map container.

\section{Accurate numerical integration of the Kolmogorov forward equation}
\label{sec:nummeth}

The PDE governing the evolution of the PDF in the present problem is given by \eqref{eq:PDEgovPDF}.
If $Q=0$, the equation is hyperbolic; if $Q>0$, the equation, strictly speaking, changes type to elliptic.
In practice, however, $Q$ is usually relatively small.
It is thus fitting to design a numerical method for accurate simulation of \eqref{eq:PDEgovPDF} based on a proven algorithm for accurate simulation of hyperbolic PDEs.
Fortunately, the fluid mechanics community has focused on the development of high performance computing techniques for numerical simulation of such ``convection-dominated''
problems for over 40 years, and these techniques are now quite refined and well understood.
The numerical method best suited to the present problem is somewhat involved; a comprehensive review of this class of methods is given in LeVeque (2002).
To focus this discussion, consider first the two-dimensional, linear, hyperbolic PDE in conservation form
\begin{equation}
\Dpartial{p(x,y,t)}{t} = - \Dpartial{\,u(x,y)\,p(x,y,t)}{x} - \Dpartial{\,v(x,y)\,p(x,y,t)}{y},
\label{eq:hyper_con}
\end{equation}
noticing that higher-dimensional cases follow as an obvious extension.
Following Chapters 4, 6, 9, 19, and 20 of LeVeque (2002), we implement a Godunov-type finite volume method by writing \eqref{eq:hyper_con} on a uniform Cartesian 2D mesh
(with constant $\Delta x$ and $\Delta y$) in the form
\begin{equation}
\frac{p^{n+1}_{ij}-p^n_{ij}}{\Delta t} = - \frac{F^n_{i+1/2,j} - F^n_{i-1/2,j}}{\Delta x} - \frac{G^n_{i,j+1/2} - G^n_{i,j-1/2}}{\Delta y},
\label{eq:hyper_conD}
\end{equation}
where the fluxes $F^n_{i-1/2,j}$ and $G^n_{i,j-1/2}$ are determined, for all $i$ and $j$, by first initializing
\begin{align*}
F^n_{i-1/2,j} &= u^+_{i-1/2,j} p^n_{i-1,j} + u^-_{i-1/2,j} p^n_{i,j},\\
G^n_{i,j-1/2} &= v^+_{i,j-1/2} p^n_{i,j-1} + v^-_{i,j-1/2} p^n_{i,j},
\end{align*}
where $u^+ = \max(u,0)$,\ \  $u^- = \min(u,0)$, etc, then applying the
{\it corner transport upwind} (CTU) terms by updating, for all $i$ and $j$,
\begin{alignat*}{2}
F^n_{i-1/2,j-1}  &\leftarrow  F^n_{i-1/2,j-1} &&-     \Delta t \frac{u^-_{i-1/2,j-1} \,v^-_{i,j-1/2}}{2} \,\frac{\Delta p^n_{i,j-1/2}}{\Delta y}, \\
F^n_{i+1/2,j-1}  &\leftarrow  F^n_{i+1/2,j-1} &&-     \Delta t \frac{u^+_{i+1/2,j-1} \,v^-_{i,j-1/2}}{2} \,\frac{\Delta p^n_{i,j-1/2}}{\Delta y}, \\
F^n_{i-1/2,j}  \ &\leftarrow\ F^n_{i-1/2,j}   &&- \ \ \Delta t \frac{u^-_{i-1/2,j}   \,v^+_{i,j-1/2}}{2} \,\frac{\Delta p^n_{i,j-1/2}}{\Delta y}, \\
F^n_{i+1/2,j}  \ &\leftarrow\ F^n_{i+1/2,j}   &&- \ \ \Delta t \frac{u^+_{i+1/2,j}   \,v^+_{i,j-1/2}}{2} \,\frac{\Delta p^n_{i,j-1/2}}{\Delta y}, \\
G^n_{i-1,j-1/2}  &\leftarrow  G^n_{i-1,j-1/2} &&-     \Delta t \frac{v^-_{i-1,j-1/2} \,u^-_{i-1/2,j}}{2} \,\frac{\Delta p^n_{i-1/2,j}}{\Delta x}, \\
G^n_{i-1,j+1/2}  &\leftarrow  G^n_{i-1,j+1/2} &&-     \Delta t \frac{v^+_{i-1,j+1/2} \,u^-_{i-1/2,j}}{2} \,\frac{\Delta p^n_{i-1/2,j}}{\Delta x}, \\
G^n_{i,j-1/2}  \ &\leftarrow\ G^n_{i,j-1/2}   &&- \ \ \Delta t \frac{v^-_{i,j-1/2}   \,u^+_{i-1/2,j}}{2} \,\frac{\Delta p^n_{i-1/2,j}}{\Delta x}, \\
G^n_{i,j+1/2}  \ &\leftarrow\ G^n_{i,j+1/2}   &&- \ \ \Delta t \frac{v^+_{i,j+1/2}   \,u^+_{i-1/2,j}}{2} \,\frac{\Delta p^n_{i-1/2,j}}{\Delta x}, 
\end{alignat*}
where $\Delta p^n_{i-1/2,j} = p^n_{ij} - p^n_{i-1,j}$,\ \  $\Delta p^n_{i,j-1/2} = p^n_{ij} - p^n_{i,j-1}$,\ \ 
then applying the {\it high-resolution correction terms} by updating, for all $i$ and $j$,
\begin{align*}
F^n_{i-1/2,j}  &\leftarrow  F^n_{i-1/2,j} + \Delta t \frac{|u_{i-1/2,j}|}{2} \left( \frac{\Delta x}{\Delta t} - |u_{i-1/2,j}|\right)\cdot\\
               & \hskip0.9in \frac{\Delta p^n_{i-1/2,j}}{\Delta x} \,\phi(\theta^n_{i-1/2,j}), \\
G^n_{i,j-1/2}  &\leftarrow  G^n_{i,j-1/2} + \Delta t \frac{|v_{i,j-1/2}|}{2} \left( \frac{\Delta y}{\Delta t} - |v_{i,j-1/2}|\right)\cdot\\
			   & \hskip0.9in \frac{\Delta p^n_{i,j-1/2}}{\Delta y} \,\phi(\theta^n_{i,j-1/2}), 
\end{align*}
where
\begin{align*}
\theta^n_{i-1/2,j}&=\begin{cases} {\Delta p^n_{i-3/2,j}}/{\Delta p^n_{i-1/2,j}} \quad &\textrm{if}\ u_{i-1/2,j}\ge 0, \\
                                  {\Delta p^n_{i+1/2,j}}/{\Delta p^n_{i-1/2,j}} \quad &\textrm{if}\ u_{i-1/2,j}<0,       \end{cases}\\
\theta^n_{i,j-1/2}&=\begin{cases} {\Delta p^n_{i,j-3/2}}/{\Delta p^n_{i,j-1/2}} \quad &\textrm{if}\ v_{i,j-1/2}\ge 0, \\
								  {\Delta p^n_{i,j+1/2}}/{\Delta p^n_{i,j-1/2}} \quad &\textrm{if}\ v_{i,j-1/2}<0,       \end{cases}
\end{align*}
and the {\it flux limiter} function $\phi(\theta)\in[0,2]$ is selected as one of several possible choices, including the {\it monotonized central-difference} ({\it MC})
limiter and the {\it van Leer} limiter:
\begin{alignat*}{2}
&\textrm{{\it MC}:}\quad &&\phi(\theta)=\max\{0,\min[(1+\theta)/2,2,2\theta]\},  \\
&\textrm{{\it van Leer}:}\quad &&\phi(\theta)=(\theta + |\theta|)/(1+|\theta|).       
\end{alignat*}
Note that exact conservation of the discrete approximation of the integral of $p$ over phase space,
as implied by the continuous formulation in \eqref{eq:hyper_con}, follows immediately from \eqref{eq:hyper_conD}.

\subsection{Numerical analysis}
\label{sec:accuracy}

In regions characterized by smooth variation of $p$, $\theta\approx 1$ and $\phi(\theta)\approx 1$, and the algorithm
described in \S \ref{sec:nummeth} is amenable to straightforward numerical analysis.  For simplicity, consider here the 1D test problem
\begin{equation}
  \Dpartial{p}{t} = - u \Dpartial{p}{x} \label{eq:analysisgod}
\end{equation}
where $u$ is a positive or negative constant.  In this case, the discretization described above reduces to
\begin{equation*}
\frac{p^{n+1}_{i}-p^n_{i}}{\Delta t} = - \frac{F^n_{i+1/2} - F^n_{i-1/2}}{\Delta x}
\end{equation*}
where
\begin{equation*}
F^n_{i-1/2}=\frac{u}{2}(p^n_i+p^n_{i-1})-\frac{u^2}{2}\frac{\Delta t}{\Delta x} (p^n_i-p^n_{i-1}),
\end{equation*}
and thus
\begin{equation*}
\frac{p^{n+1}_{i}-p^n_{i}}{\Delta t} = - u \frac{(p^n_{i+1}-p^n_{i-1})}{2\,\Delta x} + \frac{u^2\,\Delta t}{2} \frac{(p^n_{i+1}-2p^n_i + p^n_{i-1})}{(\Delta x)^2}.
\end{equation*}
Now applying to this equation the multidimensional Taylor series expansion,
\begin{align*}
y^{n+m}_{i+k} &= y^{n}_{i} + m \Delta t \big(\pderiv{y}{t}\big)^{n}_{i} + k \Delta x \big(\pderiv{y}{x}\big)^{n}_{i} + \frac{(m \Delta t)^{2}}{2} \big(\pderivsq{y}{t}\big)^{n}_{i} \\
	          &\quad + \frac{(k \Delta x)^{2}}{2} \big(\pderivsq{y}{x}\big)^{n}_{i} + m \Delta t \,\,k \Delta x \big(\frac{\partial^2 y}{\partial x \partial t}\big)^{n}_{i} + ...,
\end{align*}
and rearranging appropriately, gives
\begin{align*}
(p_t)^n_i &= - u (p_x)^n_i -  \frac{\Delta t}{2} (p_{tt})^n_i + \frac{u^2 \Delta t}{2} (p_{xx})^n_i\\
&\quad + {\mathcal O}((\Delta t)^2,(\Delta x)^2,\Delta x\,\Delta t).
\end{align*}
Differentiating \eqref{eq:analysisgod} with respect to $t$ and inserting \eqref{eq:analysisgod} into the RHS of the result, it is seen that the second and third terms on the RHS of the above expression cancel.
Thus, in regions of smooth variation of $p$, the proposed scheme is {\it second-order accurate in both space and time}.\footnote{Meaning that the error is bounded by a term proportional to $(\Delta x) ^2$ in space and $(\Delta t)^2$ in time, giving convergence of $\mathcal{O}(1/N^2)$.}
A similar analysis follows for problems in higher dimensions.

\subsection{Accounting for diffusion}
\label{sec:nummethdiff}

A diffusion term is easily added to the discretization given in \eqref{eq:hyper_conD} in a second-order central finite difference fashion simply by updating the fluxes such that, for all $i$ and $j$,
\begin{align*}
F^n_{i-1/2,j}  &\leftarrow  F^n_{i-1/2,j} + \bar\mu \frac{\Delta p^n_{i-1/2,j}}{\Delta x}, \\
G^n_{i,j-1/2}  &\leftarrow  G^n_{i,j-1/2} + \bar\mu \frac{\Delta p^n_{i,j-1/2}}{\Delta y},
\end{align*}
where $\bar\mu$ is the coefficient of the diffusion term that is applied numerically.
The flux limiter functions mentioned at the end of \S \ref{sec:nummeth} are designed to reduce the algorithm, locally, to a first-order spatial behavior while applying sufficient numerical diffusion in regions of large local curvature of $p$ on the grid, to provide a total variation diminishing (TVD) solution (that is, preventing spurious oscillations with new local minima and maxima).
We may compensate for the diffusion introduced by the numerical discretization of the convective terms simply by appropriately reducing the diffusion $\bar{\mu}$ applied in the numerical simulation of \eqref{eq:PDEgovPDF}.

\subsection{Validation}
\begin{figure}
\centering
	\includegraphics[width=0.45\textwidth]{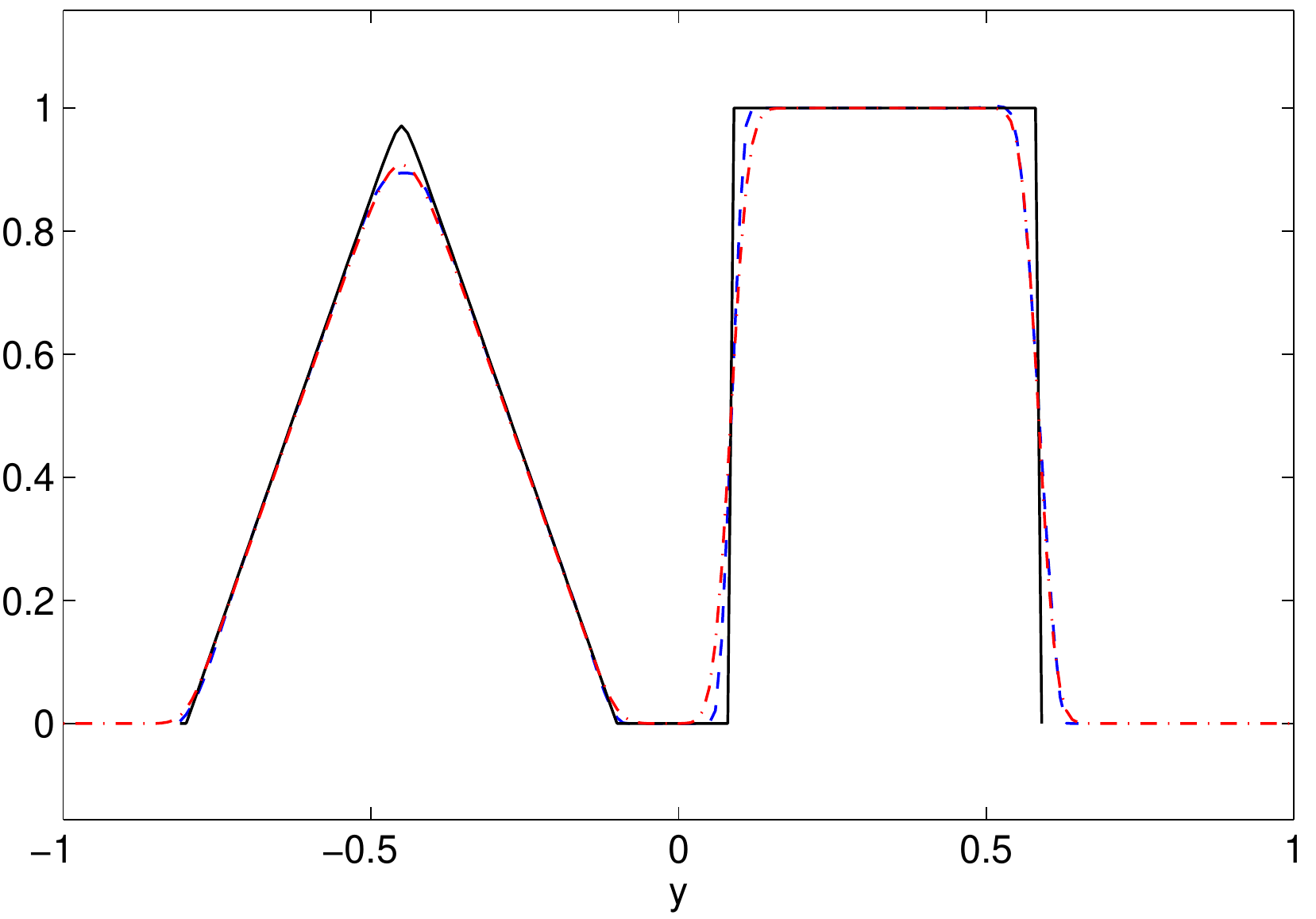}
\caption{
Computation of a sensitive test system (governed by \eqref{eq:hyper_con} with $u(x,y)=y$ and $v(x,y)=-x$), using the method in \S \ref{sec:nummeth}, taking $Q=0$, $\Delta x=0.01$, $\Delta t=0.001$, and a threshold of $10^{-3}$ (that is, tracking numerically only those cells with $p>0.001$ and their immediate neighbors).
The distributions are compared in cross section in the plane $x=0$, with the black solid line denoting the exact solution at $t=2\pi$ for $Q=0$, the blue dashed line denoting the numerical solution at $t=2\pi$ using $Q=0$, $\Delta x=0.01$, $\Delta t=0.001$ and the red dot-dashed line denoting the exact solution at $t=2\pi$ for $Q=8\times 10^{-5}$.
The leading-order error in the numerical solution at $t=2\pi$ is a small amount of diffusion of $p$.}
\label{fig:sparsity}
\end{figure}

A simple yet sensitive numerical test of the algorithm 
is given in Figure \ref{fig:sparsity}; this numerical test was taken with $u=y$ and $v=-x$ in order to give simple solid body rotation about the origin, as suggested by LeVeque (2002).
If $Q=0$, the exact solution of the test problem considered in Figure \ref{fig:sparsity} at $t=2\pi$, after a single rotation of the system about the origin, is simply the initial condition.
For the case \mbox{$Q=2 \mu I$} where $\mu$ is constant and positive,
 the exact solution of this problem at $t=2\pi$ may be obtained analytically by means of Fourier transforms.
As seen by comparing Figure \ref{fig:sparsity} to Figure 20.5 of LeVeque (2002), the result obtained via the GBEES approach is essentially identical to that obtained using the complete grid when sufficiently small threshold, time step and state-space discretization is used.
The information loss due to the discretization scheme may be quantified via the Kullback-Liebler divergence (Kullback \& Liebler, 1951), $D_{KL}(P_0,P_\epsilon)$, where a distribution $P_\epsilon$ is used to approximate the true distribution $P_0$. The Kullback-Liebler divergence using the simulation in Figure \ref{fig:sparsity} when compared to the true analytic solution is $0.089$ bits (with the distributions normalised to integrate to unity), whereas the divergence from the true case to the diffusion case calculated analytically for $Q=8\times 10^{-5}$ is $0.085$ bits\footnote{This value of $Q$ was chosen so that the divergence from the true case to the diffusion case was close to that of the true case to the numerical solution.}. The divergence from the true solution to a simulation using a truncation threshold of $10^{-16}$ (not shown) is almost the same as the more aggressively truncated example and is visually indistinguishable.
As evident by comparing the numerical solution at $t=2\pi$ in Figure \ref{fig:sparsity} to the initial condition, the discretization described in \S \ref{sec:nummeth} introduces a small numerical error in regions of high curvature.
However, by comparing the numerical solution, for $\mu=0$, $\Delta y=0.001$ and $\Delta x=0.01$, to the exact solution, for $\mu=4\times 10^{-5}$, it is evident that the leading-order error of the numerical discretization is just a bit of additional diffusion, the level of which may be determined by a suitable minimisation process.

\section{Numerical results}
\begin{figure}
\centering
\includegraphics[width=0.45\textwidth]{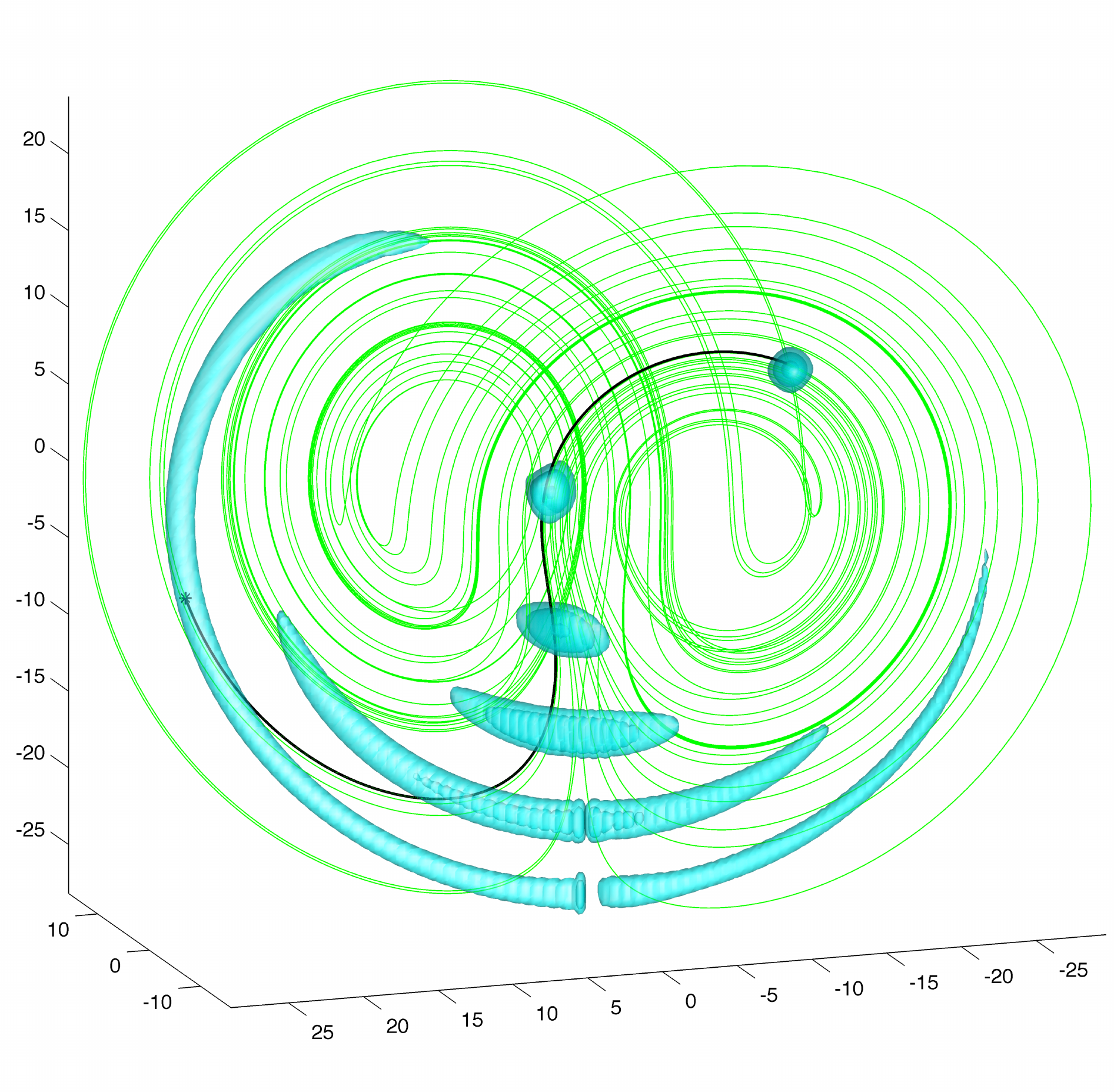}
\caption{The rapid transformation of a PDF from Gaussian to highly non-Gaussian in the Lorenz system, with no measurement updates.  Visualized are $p=0.005$, $p=0.0005$, and $p=0.00005$ isosurfaces of the PDF in phase space at $t=0$ (in the upper-right), $t=0.2$, $t=0.4$, $t=0.6$, $t=0.8$, and $t=1$. The simulation was performed with a variable time step $\Delta t\leq 0.001$, a grid spacing of $\Delta x=0.25$, and a threshold of $\epsilon=10^{-6}$. This simulation required less than 40 seconds of computation on a 2009 vintage Apple laptop computer (2.4 GHz Intel Core 2 Duo) using a single-threaded C++ implementation of the present algorithm tracking about 50,000 active cells (of $5.3\times 10^6$ total cells in the domain shown) at $t=1$. The modest memory requirements are proportional to the number of cells.
}
\label{fig:transformation}
\end{figure}


\begin{figure}
\centering
\includegraphics[width=0.45\textwidth]{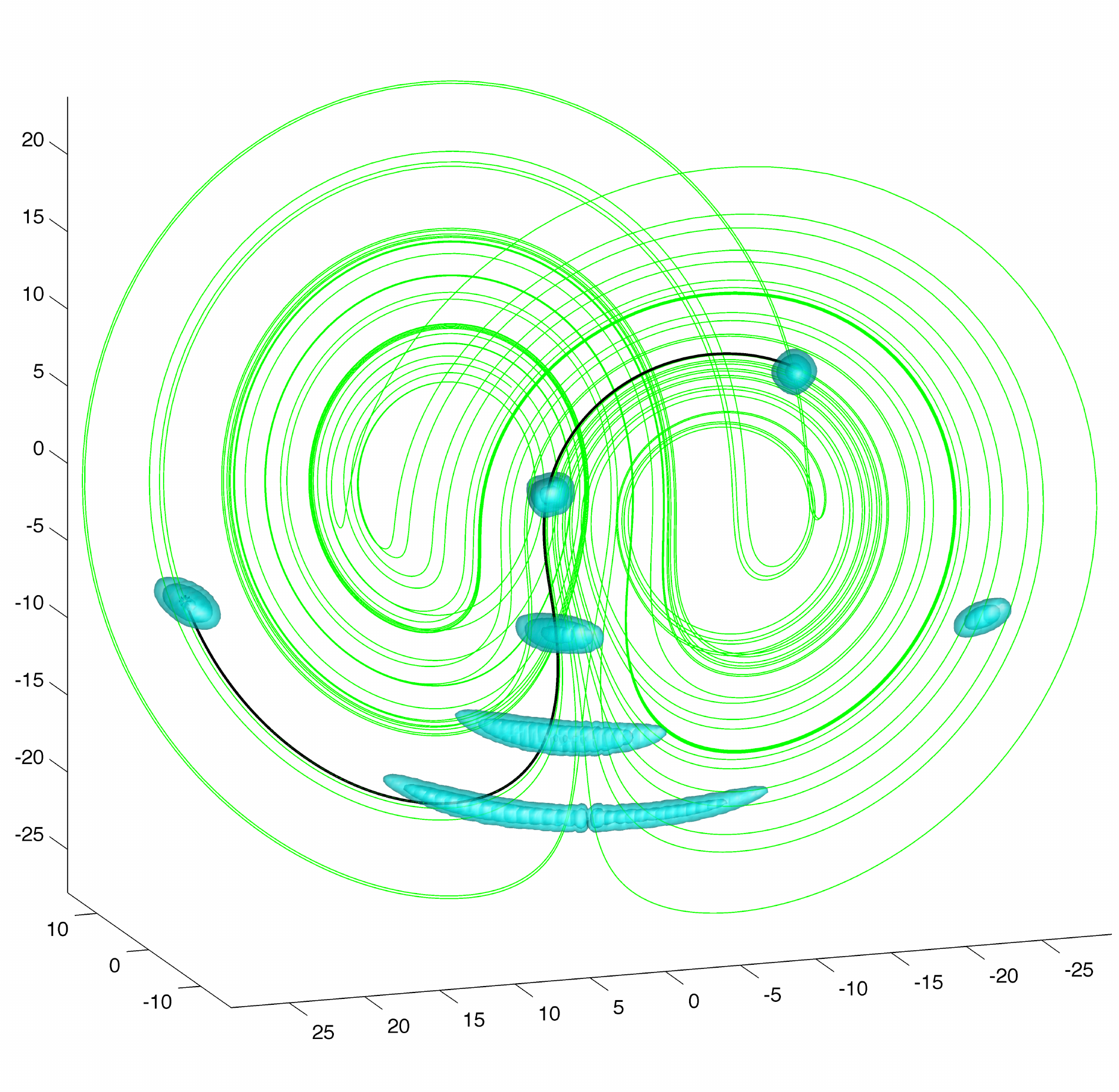}
\caption{The same simulation as Figure \ref{fig:transformation} but with measurements of $x_3$ (vertical axis), with Gaussian uncertainty, at every time step. The black line is the `true' state which generates the measurements. Note that, by $t=1$, the PDF splits into two concentrated regions. This simulation required about 4 seconds of computation with the same hardware, tracking about 4,000 active cells (of $5.3\times 10^6$ total cells in the domain shown) at $t=1$.}\label{fig:measurements}
\end{figure}

A Bayesian approach is justified when the uncertainty of the estimate is significantly non-Gaussian, such as in the estimation of a nonlinear system with relatively large uncertainty, leaving us with particle filtering or grid-based methods; what is perhaps still uncertain is the numerical tractability of a grid-based approach when one exploits the sparsity of the PDF in the manner described in \S \ref{sec:nummethsparse}.
Thus, in order to test the efficiency of the GBEES algorithm, as well as to demonstrate how it can capture with unprecedented accuracy the evolution of a non-Gaussian PDF, we have applied the GBEES algorithm to the estimation of the three-state Lorenz system
\begin{equation*}
\frac{d\x}{dt}=\f(\x), \quad
\x = \begin{pmatrix} {x}_1 \\ {x}_2 \\ {x}_3 \end{pmatrix}, \quad
\f(\x) = \begin{pmatrix} \sigma\,( x_2 - x_1 )\\ -   x_2 - x_1 x_3\\ -b\,x_3 + x_1 x_2 - b\,r \end{pmatrix}
\end{equation*}
with $\sigma=4$, $b=1$, and $r=48$.  For these parameter values, the system is chaotic, and the attractor takes the familiar form indicated by the green line in Figure \ref{fig:transformation}.  Also illustrated in this Figure \ref{fig:transformation} is the evolution of an initially (at $t=0$) Gaussian PDF $p_{\x}(\x',t)$, the evolution of which is governed by the Kolmogorov equation \eqref{eq:PDEgovPDF}, with no measurement updates applied and no added process noise (diffusion). The distribution narrows significantly in the direction normal to the attractor, and spreads out rapidly in the direction of the maximum local Lyapunov exponent along the attractor; by $t=1$, the PDF is highly non-Gaussian.
Note also in the $t=0.8$ and $t=1$ isosurfaces the remarkable division of the PDF into two distinct lobes in the vicinity of the $x_3$ axis (the vertical coordinate axis in the figures), which is invariant and unstable in the Lorenz system.

Figure \ref{fig:measurements} represents the evolution of the PDF when measurements (with Gaussian uncertainty) of $x_3$ are taken at every time step. Computationally, the problem addressed in the figure is significantly easier than the ``open-loop'' problem illustrated in Figure \ref{fig:transformation}, as the number of active cells by $t=1$ is reduced from 50,000 to only 4,000; the computation time for this simulation is also reduced accordingly, from 40 to 4 seconds for the time interval shown.
The PDF at time $t=1$ splits into two concentrated regions on the left and right sides of the figure.
Future measurements might reveal in which region the state really is; until such measurements are received, the GBEES algorithm is perfectly capable of following both. A plain Kalman filter, which assumes a central estimate, would necessarily fail to model such a splitting.

\section{Analysis and Conclusions}


A novel algorithm is introduced in this paper to exploit the remarkable sparsity of the evolving PDF in phase space, thereby, for the first time, making high-resolution grid-based Bayesian estimation computationally tractable for nontrivial systems. The method generalises straightforwardly to any number dimensions, with computational cost expected to be a trade-off between the curse of dimensionality and the increased sparseness of the PDF.
In application, the algorithm developed is shown to track, with unprecedented fidelity, the completely non-Gaussian PDF of the estimate of a Lorenz system, both with and without measurement updates.
The simulation exhibits a competition between information loss due to the random state disturbances and stretching of the PDF in the unstable directions of the system, and information gain from measurements.

Grid-based Bayesian estimation algorithms are sometimes referred to as \emph{approximate} grid-based methods. We point out that the numerical analysis of \S \ref{sec:accuracy} establishes that the numerical method used to propagate the Kolmogorov equation in the present grid-based estimation algorithm is \emph{second-order accurate in both space and time}; this compares favorably to the {(\it slower than linear)} ${\mathcal O}(1/\sqrt{N})$ convergence rate of particle methods applied to the Kolmogorov equation (see Bernard, Talay, \& Tubaro 1994).

Finally, Lagrangian (that is, particle-based) simulation techniques have been explored for decades in the field of fluid mechanics, but for $n>2$ remain mostly a research novelty.
On the other hand, grid-based methods (often with adaptive grids to focus the computational effort where it is needed) have proven immensely successful in a variety of complex situations in fluid mechanics, such as in the characterization of fluid turbulence and in the design of commercial airliners, where computational methods have largely supplanted repetitive wind-tunnel testing.
There appears to be no reason why the same success of grid-based methods will not also be realized in Bayesian estimation approaches, once the community working on such problems fully appreciate how the remarkable sparsity of the PDF in such problems may be exploited.

\subsection*{Acknowledgements}

The authors gratefully acknowledge Prof.~Paulo Luchini for insightful discussions related to this work.
An Imperial College Junior Research Fellowship (AS) is also gratefully acknowledged.

\end{document}